\documentclass{article}
\usepackage{amssymb}
\usepackage{amsfonts}
\usepackage{amsmath}
\usepackage{mathptmx}
\usepackage{geometry}
\usepackage{graphicx}
\usepackage{latexsym}
\usepackage{hyperref}

\newtheorem{theorem}{Theorem}
\newtheorem{acknowledgement}{Acknowledgement}
\newtheorem{algorithm}{Algorithm}

\newtheorem{definition}{Definition}

\newtheorem{proposition}{Proposition}
\newtheorem{remark}{Remark}

\newenvironment{proof}[1][Proof]{\noindent\textbf{#1.} }{\ \rule{0.5em}{0.5em}}

\begin{document}

\title{Characterization of Positive Operators}
\author{L\'{u}cio Fassarella \\
{\small DMA/CEUNES/UFES}}
\maketitle

\begin{abstract}
A characterization of positive operators on finite dimensional complex
vector spaces based on the Routh-Hurwitz criterion.
\end{abstract}


\section{Introduction}

Here we give a characterization of the positive operators within the class of self-adjoint operators on finite dimensional complex vector spaces. The result is a direct consequence of the \textit{Routh-Hurwitz theorem} \cite{hurwitz} which establishes a sufficient condition  for all eigenvalues of a real polynomial to have negative real part.\\

In addition to providing a test for the positivity of self-adjunct operators, Theorem \ref{teorema_caracterizacao-op} can be applied to problems in which positive operators occur as unknown variables. The need to deal with this situation appeared in a research whose main objective was to find sufficient conditions for the conservation of energy in dissipative quantum systems \cite{lucio2011}.\\

About the structure of this article, Section \ref{secao_op} reviews basic facts about positive operators and characteristic polynomials; Section \ref{secao_rh} presents the \textit{Routh-Hurwitz criterion} and an analogous result which determines whether all (complex) roots of a real polynomial have non-negative real part; in the Section \ref{secao_tcop}, our main theorem is stated and demonstrated, combining facts presented earlier; finally, Section \ref{secao_exemplos} presents the application of the main theorem in two and three dimensions.

\section{Positive operators\label{secao_op}}

Let $\left( \mathcal{H},\left\langle ,\right\rangle \right) $ be a Hilbert
space, where $\left\langle ,\right\rangle $ is a inner product (anti-linear
in the first entry and linear in the second entry). $\left\Vert \
\right\Vert $ denotes the norm associated to this inner product:%
\begin{equation*}
\left\Vert \varphi \right\Vert =\sqrt{\left\langle \varphi ,\varphi
\right\rangle },\ \forall \varphi \in \mathcal{H}.
\end{equation*}%
Here, an `operator in $\mathcal{H}$' means a bounded linear operator from $%
\mathcal{H}$ to $\mathcal{H}$. In definitions below, I follow the references \cite{Halmos}, \cite{Halmos1}
and \cite{Conway}. (After some generalities, $\mathcal{H}$ will
mean $\mathbb{C}^{n}$, with $n$ being a positive integer
number.)

\begin{definition}[adjoint, self-adjoint\ operator, positive operator]
\label{def_adjoint} Let $A$ be an operator in $\mathcal{H}$.\newline
The adjoint operator of $A$ is the operator $A^{\ast }:\mathcal{H}%
\rightarrow \mathcal{H}$ such that:\footnote{%
The existence of the operator $A^{\ast }$ follows from \textit{Riesz
representation theorem} \cite[p.39]{Halmos1} \cite[p.12, p.31]{Conway}.}%
\begin{equation*}
\left\langle A^{\ast }\varphi ,\psi \right\rangle =\left\langle \varphi
,A\psi \right\rangle ,\ \forall \varphi ,\psi \in \mathcal{H}.
\end{equation*}%
$A$ is self-adjoint when $A^{\ast }=A$, \textit{i.e.},%
\begin{equation*}
\left\langle \varphi ,A\psi \right\rangle =\left\langle A\varphi ,\psi
\right\rangle ,\ \forall \varphi ,\psi \in \mathcal{H}.
\end{equation*}%
$A$ is positive when%
\begin{equation*}
\left\langle \varphi ,A\varphi \right\rangle \geq 0 ,\ \forall \varphi
\in \mathcal{H}.
\end{equation*}
\end{definition}

\begin{proposition}
Any positive operator is self-adjoint also.
\end{proposition}

This proposition follows from \textit{polarization identity} \cite[p.171-172]%
{Halmos} \cite[p.4]{Conway}, which provides the inner product in terms of the
norm:
\begin{equation*}
\left\langle \psi ,\varphi \right\rangle =\frac{1}{4}\left\{ \left(
\left\Vert \psi +\varphi \right\Vert ^{2}-\left\Vert \psi -\varphi
\right\Vert ^{2}\right) -i\left( \left\Vert \psi +i\varphi \right\Vert
^{2}-\left\Vert \psi -i\varphi \right\Vert ^{2}\right) \right\} ,\
\forall \psi ,\varphi \in \mathcal{H}.
\end{equation*}

Our goal is to find a criterion for characterizing positive operators within the set of self-adjoint operators. The argument of Section \ref{secao_tcop} uses some definitions and results of \textit{Spectral Theory}, the relevant facts of which are recalled below.

\medskip

\begin{definition}[eigenvalue, eigenvector, characteristic polynomial]
Let $A$ be an operator in $\mathcal{H}$. An `eigenvalue of $A$' is a number $%
\lambda \in \mathbb{C}$ such that $\ker \left( \lambda I-A\right) \neq
\left\{ 0\right\} $; in this case, the elements of $\ker \left( \lambda
I-A\right) \setminus \left\{ 0\right\} $ are called `eigenvectors of $A$
associated to the eigenvalue $\lambda $'. Explicitly: $\varphi \in \mathcal{H%
}$ is an eigenvector of $A$ associated to the eigenvalue $\lambda \in 
\mathbb{C}$ if, and only if, 
\begin{equation*}
\varphi \neq 0\ \ \text{and}\ \ A\varphi =\lambda \varphi.
\end{equation*}%
The `spectrum of $A$' is the set of its eigenvalues:%
\begin{equation*}
\sigma _{A}:=\left\{ \text{eigenvalues of }A\right\}.
\end{equation*}%
When $\mathcal{H}$ has finite dimension, the `characteristic polynomial of $%
A $' is defined by%
\begin{equation}
p_{A}\left( z\right) :=\det \left( zI-A\right).
\label{form_polinomio-characteristic}
\end{equation}%
In this situation, the spectrum of $A$ is exactly the set of roots of the
characteristic polynomial of $A$:%
\begin{equation*}
\sigma _{A}=\left\{ z \in \mathbb{C} / p_{A}\left( z\right)= 0 \right\}.
\end{equation*}
\end{definition}

\begin{proposition}
\label{teorema_autovalores} Let $A$ be an operator in the Hilbert space $%
\mathcal{H}$.

i) If $A$ is self-adjoint, then its eigenvalues are real numbers;

ii)\ If $A$ is positive, then its eigenvalues are non-negative real numbers.

\begin{proof}
Let $\lambda $ be an eigenvalue of $A$ and let $\varphi \in \mathcal{H}$ ($%
\varphi \neq 0$) an eigenvector of $A$ associated to $\lambda $.\newline
If $A$ is self-adjoint, then 
\begin{equation*}
\lambda =\frac{\left\langle \varphi ,\lambda \varphi \right\rangle }{%
\left\Vert \varphi \right\Vert ^{2}}=\frac{\left\langle \varphi ,A\varphi
\right\rangle }{\left\Vert \varphi \right\Vert ^{2}}=\frac{\left\langle
A\varphi ,\varphi \right\rangle }{\left\Vert \varphi \right\Vert ^{2}}=\frac{%
\left\langle \lambda \varphi ,\varphi \right\rangle }{\left\Vert \varphi
\right\Vert ^{2}}=\bar{\lambda}.
\end{equation*}%
If $A$ is positive, then%
\begin{equation*}
\lambda =\frac{\left\langle \varphi ,\lambda \varphi \right\rangle }{%
\left\Vert \varphi \right\Vert ^{2}}=\frac{\left\langle \varphi ,A\varphi
\right\rangle }{\left\Vert \varphi \right\Vert ^{2}}\geq 0.
\end{equation*}
\end{proof}
\end{proposition}

\begin{theorem}[Spectral]
\label{teorema_espectral} Assume that $\mathcal{H}$ is separable and let $A$
be a (bounded) operator in $\mathcal{H}$. Then, $A$ is self-adjoint if, and only if, there exists a sequence of
orthogonal projections $\left( P_{1},P_{2},...\right) $%
\begin{equation*}
P_{k}P_{j}=P_{k}P_{j}=\delta _{kj}P_{k} ,\ \forall k\neq j=1,2,...
\end{equation*}%
and a sequence of real numbers $\ \left(
\lambda _{1},\lambda _{2},...\right) $ such that\footnote{%
When the dimension of $\mathcal{H}$ is infinite, the sum can be a convergent serie (w.r.t. the norm topology).}%
\begin{equation*}
I=\sum_{k}P_{k} ,\ A=\sum_{k}\lambda _{k}P_{k}.
\end{equation*}%
In this case%
\begin{equation*}
\sigma _{A}=\left\{ \lambda _{1},\lambda _{2},...\right\}.
\end{equation*}%
and%
\begin{equation*}
\emph{Im}(P_{k})=\ker \left( \lambda _{k}I-A\right) \ \ ,\ \forall k=1,2,...
\end{equation*}
\end{theorem}

\begin{proposition}\label{operadores-positivos}
\label{teorema_operadores-positivos-polinomio-caracteristico}\ Assume $%
\mathcal{H}$ is separable and let $A$ be a self-adjoint operator in $%
\mathcal{H}$. Then, $A$ is positive if, and only if, its eigenvalues are non-negatives.

\begin{proof}
The implication was already demonstrated in the Proposition \ref%
{teorema_autovalores}; so, here I prove the reciprocal sentence (even though
it also follows from the identity below). Take the spectral decomposition of 
$A$ according with Theorem \ref{teorema_espectral}; then, for all $\varphi
\in \mathcal{H}$ it holds%
\begin{equation*}
\left\langle \varphi ,A\varphi \right\rangle =\left\langle \varphi
,\sum_{k}\lambda _{k}P_{k}\varphi \right\rangle =\sum_{k}\lambda
_{k}\left\langle \varphi ,P_{k}\varphi \right\rangle =\sum_{k}\lambda
_{k}\left\langle P_{k}\varphi ,P_{k}\varphi \right\rangle =\sum_{k}\lambda
_{k}\underbrace{\left\Vert P_{k}\varphi \right\Vert ^{2}}_{\geq 0}.
\end{equation*}%
Therefore, if all eigenvalues of $A$ are non-negatives, the last sum above
is non-negative also. This proves that $A$ is positive, according with
Definition \ref{def_adjoint}.
\end{proof}
\end{proposition}

\begin{proposition}
	\label{teorema_polinomio-caracteristico} Assume that $\mathcal{H}$ is finite dimensional and let $A$ be a self-adjoint operator on $\mathcal{H}$. Then, the characteristic polynomial of $A$ is real (\textit{i.e.}, its has real coefficients).
	
	\begin{proof}
		Let $\bar{p}_{A}\left( z\right) $ be the polynomial whose coefficients are 		the corresponding complex conjugated of the coefficients of the
		characteristic polynomial $p_{A}\left( z\right) $; then:%
		\begin{equation*}
		\bar{p}_{A}\left( z \right) = \overline{p_{A}\left( \bar{z} \right) } ,\
		\forall z\in \mathbb{C}.
		\end{equation*}%
		Using the algebraic properties of the conjugation $\ast $ (which assigns
		each operator to its adjoin), it follows:%
		\begin{equation*}
		\overline{p_{A}\left( z\right) }=\overline{\det \left( zI-A\right) }=\det
		\left( zI-A\right) ^{\ast }=\det \left( \bar{z}I-A\right) =p_{A}\left( \bar{z%
		}\right) ,\ \forall z\in \mathbb{C}.
		\end{equation*}%
		Combining the two identities above, it follows:%
		\begin{equation*}
		\bar{p}_{A}\left( x\right) =p_{A}\left( x\right) ,\ \forall x\in \mathbb{R}.
		\end{equation*}%
		Therefore, the polynomials $\bar{p}_{A}\left( z\right) $ and $p_{A}\left(
		z\right) $ have the same restriction to $\mathbb{R}$; this implies that they
		have the same coefficients, from what follows that the coefficients of $%
		p_{A}\left( z\right) $ are real numbers.
	\end{proof}
\end{proposition}

From Theorem \ref{operadores-positivos} and Proposition \ref{teorema_polinomio-caracteristico}, we see that one can obtain a characterization of the positive operators within the class of self-adjoint operators in $\mathbb{C} ^{n}$ using any characterization of the real polynomials whose complex roots have non-negative real part. This points to the \textit{Routh-Hurwitz
criterion}.

\section{Routh-Hurwitz criterion\label{secao_rh}}

Let $p\left( z\right)$ be a polynomial with complex coefficients,

\begin{equation}
p\left( z\right) =b_{0}z^{n}+b_{1}z^{n-1}+...+b_{n-1}z+b_{n}.
\label{form_polinomio}
\end{equation}%
The \textit{Hurwitz determinants} of $p\left( z\right) $ are defined by%
\begin{equation}
\Delta _{k}:=\det \left( 
\begin{array}{ccccccc}
b_{1} & b_{3} & b_{5} & b_{7} & \cdots & \cdots & b_{2k-1} \\ 
b_{0} & b_{2} & b_{4} & b_{6} & \cdots & \cdots & b_{2k-2} \\ 
0 & b_{1} & b_{3} & b_{5} & \cdots & \cdots & b_{2k-3} \\ 
0 & b_{0} & b_{2} & b_{4} & \cdots & \cdots & \vdots \\ 
\vdots & 0 & b_{1} & b_{3} & b_{5} & \cdots & b_{k+2} \\ 
\vdots & \vdots & \vdots & \vdots & \vdots & \ddots & b_{k+1} \\ 
0 & 0 & \cdots & \cdots & b_{k-4} & b_{k-2} & b_{k}%
\end{array}%
\right) \ \left(k=1,2,...,n \right) \label{form_determinantes-hurwitz},
\end{equation}%
where
\begin{equation*}
b_{j}=0\ \ ,\ \forall j\in \mathbb{N},\ j>n.
\end{equation*}%
In particular, for $n=3$ this definition means%
\begin{equation*}
\Delta _{1}=b_{1}\ \ ,\ \ \Delta _{2}=\det \left( 
\begin{array}{cc}
b_{1} & b_{3} \\ 
b_{0} & b_{2}%
\end{array}%
\right) \ \ ,\ \ \Delta _{3}=\det \left( 
\begin{array}{ccc}
b_{1} & b_{3} & 0 \\ 
b_{0} & b_{2} & 0 \\ 
0 & b_{1} & b_{3}%
\end{array}%
\right).
\end{equation*}

\begin{theorem}[Routh-Hurwitz criterion]
\label{teorema_rhcriterion}\ Consider the polynomial (\ref{form_polinomio})
has real coefficients and that its leading coefficient is positive ($b_{0}>0$%
). Then, the roots of this polynomial have negative real part if, and only
if, the Hurwitz determinants of the polynomial are positive:%
\begin{equation*}
\Delta _{k}>0\ \ ,\ \forall k=1,...,n
\end{equation*}
\end{theorem}

For a proof of this theorem I refer to \cite[p.231]{gantmacher}.\\

In the sequel, I will state and demonstrate theorems analogous to the
Routh-Hurwitz criterion. By convention, the phrase \textquotedblleft $w\in \mathbb{C}$ is a root of
the polynomial $p\left( z\right) $ with multiplicity zero\textquotedblright\ means that \textquotedblleft $w$ is not a root of $p\left( z\right) $ \textquotedblright .

\begin{proposition}[Extended Routh-Hurwitz criterion]
\label{teorema_routh-hurwitz-extendido} Consider a polynomial $p\left( z\right)$ having real coefficients and positive leading coefficient. Let $n_{0}$ be the multiplicity of number zero as a root of $p\left(z\right) $, $0\leq n_{0}<n$. Then, the non-zero roots of $p\left( z\right) $ have negative real part if, and only if, its Hurwitz determinants with indexes from $1$ to $n-n_{0}$ are positive%
\begin{equation*}
\Delta _{k}>0\ \ ,\ \forall k=1,...,n-n_{0}.
\end{equation*}

\begin{proof}
Let $p\left( z\right)$ be as in Eq.\ref%
{form_polinomio}. The case $n_{0}=0$ reduces to the Routh-Hurwitz criterion (Theorem \ref%
{teorema_rhcriterion}); therefore, I will assume that $n_{0}\geq 1$. From the hypothesis, the coefficients of the monomials of $p\left( z\right) $
having degree less than $n_{0}$ are equal to zero 
\begin{equation}
b_{n-n_{0}+1}=...=b_{n}=0  \label{equac_coeficientes-nulos}.
\end{equation}%
I define the auxiliary polynomial%
\begin{equation*}
\hat{p}\left( z\right)
:=b_{0}z^{n-n_{0}}+b_{1}z^{n-n_{0}-1}+...+b_{n-n_{0}-1}z+b_{n-n_{0}}.
\end{equation*}%
Therefore 
\begin{equation*}
p\left( z\right) =z^{n_{0}}\hat{p}\left( z\right).
\end{equation*}%
The polynomial $\hat{p}\left( z\right) $ has degree $n-n_{0}$ and its roots are equal to the non-zero roots of $p\left( z\right) $; further, the Hurwitz determinants of $\hat{p}\left( z\right) $ and $p\left( z\right) $ are also equal due to Eq.\ref{equac_coeficientes-nulos} and Definition \ref%
{form_determinantes-hurwitz}. As happens with $p\left( z\right) $, $\hat{p} \left( z\right) $ is a real polynomial with positive leading coefficient; therefore, the thesis follows from Theorem \ref{teorema_rhcriterion} applied to $\hat{p}\left( z\right) $.
\end{proof}
\end{proposition}

\begin{remark}
In the case $p\left( z\right) =z^{n}$ (the number zero is a root of the
polynomial with multiplicity equal to the polynomial's degree), all Hurwitz coefficients of $p\left( z\right) $ are zero and (evidently) the polynomial has no non-zero root.
\end{remark}

\begin{proposition}[Symmetric of the extended Routh-Hurwitz criterion]
\label{teorema_routh-hurwitz-simetrico} Consider a polynomial $p\left( z\right)$ having real coefficients and positive leading coefficient. Let $n_{0}$ be the multiplicity of the number zero as a root of $p\left(z\right) $, $0\leq n_{0}<n$. Then, the non-zero
roots of $p\left( z\right) $ have positive real part if, and only if, the
Hurwitz determinants satisfy the condition:%
\begin{equation}
\left( -1\right) ^{1+\left\lfloor k/2\right\rfloor }\Delta _{k}>0\ \
\forall k=1,2,...,n-n_{0}  \label{form_serh-criterion},
\end{equation}%
where $\left\lfloor k/2\right\rfloor $ is the greatest integer less then or equal to $k/2$:%
\begin{equation*}
\left\lfloor x\right\rfloor :=\max \left\{ m\in \mathbb{Z};\ m\leq x\right\}
\ \ \forall x\in \mathbb{R}.
\end{equation*}%
Therefore, the conditions given by Eq.\ref{form_serh-criterion} are necessary and
sufficient to guarantee that all roots of the polynomial $p\left( z\right) $
have non-negative real part.

\begin{proof}
Consider the auxiliary polynomial%
\begin{equation*}
\tilde{p}\left( z\right) :=\left( -1\right) ^{n}p\left( -z\right).
\end{equation*}%
Then, the non-zero roots of $p\left( z\right) $ have positive real part if,
and only if, all non-zero roots of $\tilde{p}\left( z\right) $ have negative
real part. The coefficients of $\tilde{p}\left( z\right) $ are given in
terms of the coefficients of $p\left( z\right) $ by%
\begin{equation*}
\tilde{b}_{k}=\left( -1\right) ^{k}b_{k}\ \ \forall k=1,2,...,n.
\end{equation*}%
Therefore, the Hurwitz determinants of $p\left( z\right) $ and $\tilde{p}%
\left( z\right) $ satisfy the identities\footnote{%
The signal came from the number $1+\left\lfloor k/2\right\rfloor $ of lines
of the matrix defining $\tilde{\Delta}_{k}$ which are different by a factor $%
-1$ from the corresponding lines of the matrix defining $\Delta _{k}$: by a
basic property of the determinant, each different line implies a factor $-1$
in the relation between the two determinants.}%
\begin{equation*}
\tilde{\Delta}_{k}=\left( -1\right) ^{1+\left\lfloor k/2\right\rfloor
}\Delta _{k}\ \  \forall k=1,2,...,n.
\end{equation*}%
As with $p\left( z\right) $, $\tilde{p}\left( z\right) $ is a real
polynomial with positive leading coefficient ($\tilde{b}_{0}=b_{0}$) and the
number zero is a root of it with multiplicity $n_{0}$; therefore, the thesis
follows from Proposition \ref{teorema_routh-hurwitz-extendido} applied to $\tilde{p}\left( z\right) $.
\end{proof}
\end{proposition}

\section{Characterization of positive operators\label{secao_tcop}}

\begin{theorem}[characterization of positive operators]
\label{teorema_caracterizacao-op}\ Let $A$ be a self-adjoint operator in $%
\mathbb{C}^{n}$ and let $n_{0}\in \left\{ 0,1,...,n\right\} $ be the
multiplicity of the number zero as a root of the characteristic polynomial
of $A$. Then, $A$ is a positive operator if, and only if,%
\begin{equation*}
\left( -1\right) ^{1+\left\lfloor k/2\right\rfloor }\Delta _{k}>0\ \ 
\forall k=1,2,...,n-n_{0}.
\end{equation*}

\begin{proof}
The result follows directly from the combination of Theorem \ref%
{teorema_operadores-positivos-polinomio-caracteristico} and Theorem \ref%
{teorema_routh-hurwitz-simetrico}.
\end{proof}
\end{theorem}

Naturally, this result motivates us to know expressions for the
coefficients of the characteristic polynomial of operators. So, I highlight the following results:

\begin{proposition}[Coefficients of the Characteristic Polynomials]
\label{teorema_coeficientes-pc}\ Let $A$ be a self-adjoint operator in $%
\mathbb{C}^{n}$ and let its characteristic polynomial be%
\begin{equation*}
p_{A}\left( z\right) =\det \left( zI-A\right)
=:b_{0}z^{n}+b_{1}z^{n-1}+...+b_{n-1}z+b_{n}.
\end{equation*}%
i) Formula in terms of the traces:%
\begin{equation*}
b_{0}=1;\ b_{k}=-\frac{1}{k}\left\{ b_{k-1}\mathrm{tr}\left( A\right)
+b_{n-2}\mathrm{tr}\left( A^{2}\right) +...+b_{1}\mathrm{tr}\left(
A^{k-1}\right) +b_{0}\mathrm{tr}\left( A^{k}\right) \right\} \ \ \forall
k=1,...,n.
\end{equation*}%
ii) Formula in terms of subdeterminants:%
\begin{equation*}
b_{0}=1\ \ ;\ \ b_{k}=\left( -1\right) ^{k}\sum_{j_{1}<...<j_{k}}\det \left( 
\begin{array}{cccc}
a_{j_{1}j_{1}} & a_{j_{1}j_{2}} & \cdots & a_{j_{1}j_{k}} \\ 
a_{j_{2}j_{1}} & a_{j_{2}j_{2}} & \cdots & a_{j_{2}j_{k}} \\ 
\vdots & \vdots & \ddots & \vdots \\ 
a_{j_{k}j_{1}} & a_{j_{k}j_{2}} & \cdots & a_{j_{k}j_{k}}%
\end{array}%
\right) \ \ \forall k=1,...,n.
\end{equation*}%
In particular,%
\begin{equation*}
b_{1}=-\mathrm{tr}\left( A\right) \ \ ,\ \ b_{n}=\left( -1\right) ^{n}\det
\left( A\right).
\end{equation*}
\end{proposition}

Proofs for such formulas can be found in \cite{lewin} and \cite{collings},
respectively. Another formula can be found in \cite{pennisi1987}.\\

Finally, I collect the facts in a simple algorithm:

\begin{algorithm}
\textbf{Characterization of positive operators}\label{algoritmo}\newline
Let $A$ be a operator in $\mathbb{C}^{n}$ and consider its representation in
the canonical basis (or any other orthogonal basis) of $\mathbb{C}^{n}$:%
\begin{equation*}
\left[ A\right] =\left( a_{ij}\right) _{i,j=1,...,n}.
\end{equation*}%
i) Verify if the operator is self-adjoint, \textit{i.e.}, if $\bar{a}%
_{ij}=a_{ji}\ \forall i,j=1,...,n$;\newline
ii) Calculate the coefficients of the characteristic polynomial of $A$,
using one of the formulas in Theorem \ref{teorema_coeficientes-pc};\newline
iii) Calculate the multiplicity of the number zero as a root of the
characteristic polynomial of $A$, using, for example, the formula%
\begin{equation*}
n_{0}=\min \left\{ k\in \mathbb{N}\ /\ \left. \frac{d^{k}p_{A}}{dz^{k}}%
\right\vert _{z=0}\neq 0\right\}.
\end{equation*}%
\newline
iv) Calculate the Hurwitz determinants of the characteristic polynomial of $%
A$ up to the order $n-n_{0}$, using Definition \ref%
{form_determinantes-hurwitz};\newline
v) Verify if the Hurwitz determinants calculated satisfy the condition of Theorem
\ref{teorema_caracterizacao-op}.
\end{algorithm}

The computational efficiency of this Algorithm \ref{algoritmo} could be
compared with the computational efficiency to calculate the operator's eigenvalues with sufficient precision. I let this for the reader.

\section{Special Cases\label{secao_exemplos}}

In this section we get from Theorem \ref{teorema_caracterizacao-op} explicit conditions for a self-adjoint operator to be positive on two and three dimensions. These expressions are presented in terms of determinants and traces since we apply Definition \ref{form_determinantes-hurwitz} and Proposition \ref{teorema_coeficientes-pc}.

\subsection{Dimension $n=2$}

Let $A$ be a self-adjoint operator in $\mathbb{C}^{2}$ and consider its
matrix representation with respect to the canonical basis (or with respect
to any other orthogonal basis):%
\begin{equation*}
\left[ A\right] =\left( a_{ij}\right) _{i,j=1,2}.
\end{equation*}%
Self-adjointness means%
\begin{equation*}
a_{ij}=\bar{a}_{ji}\in \mathbb{C}\ \  \forall i,j=1,2.
\end{equation*}%
The coefficients of the characteristic polynomial of $A$ are, as one can
obtain directly from (\ref{form_polinomio-characteristic}) or from Theorem %
\ref{teorema_coeficientes-pc}:%
\begin{equation*}
b_{0}=1\ ,\ b_{1}=-\mathrm{tr}\left( A\right),\ b_{2}=\det \left( A\right).
\end{equation*}%
The Hurwitz determinants (\ref{form_determinantes-hurwitz}) are:%
\begin{equation*}
\Delta _{1}=b_{1}=-\mathrm{tr}\left( A\right) , \ \Delta
_{2}=b_{1}b_{3}-b_{0}b_{2}=\det \left( A\right).
\end{equation*}%
Now, to get the necessary and sufficient conditions for $A$ be positive, we
have to consider three different situations, distinguished by the
multiplicity, denoted here by $\mu $, of the number zero as a root of the
characteristic polynomial of $A$.

\underline{Case $\mu =0$} $\left( b_{2}=\det \left( A\right) \neq 0\right) $%
: Theorem \ref{teorema_caracterizacao-op} stablishes that $A$ is positive
if, and only if,%
\begin{equation*}
\mathrm{tr}\left( A\right) >0\ \ ,\ \det \left( A\right) >0.
\end{equation*}

\underline{Case $\mu =1 $} $\left( b_{2}=\det \left( A\right) =0\ ,\ b_{1}=-%
\mathrm{tr}\left( A\right) \neq 0\right) $: Theorem \ref%
{teorema_caracterizacao-op} stablishes that $A$ is positive if, and only if,%
\begin{equation*}
\mathrm{tr}\left( A\right) >0.
\end{equation*}

\underline{Case $\mu =2 $} $\left( b_{2}=\det \left( A\right) =0\ ,\ b_{1}=-%
\mathrm{tr}\left( A\right) =0\right) $: Theorem \ref%
{teorema_caracterizacao-op} stablishes that $A$ is positive (without any
further condition).\\

Finally, we can collect the cases in a single sentence:

\begin{quote}
A self-adjoint operator $A$ in $\mathbb{C}^{2}$ is positive if, and only if,
it satisfies one out of the two following conditions:%
\begin{equation}
\left\{ 
\begin{array}{c}
\left( i\right) \ \det \left( A\right) =0\ ,\ \mathrm{tr}\left( A\right)
\geq 0; \\ 
\left( ii\right) \ \det \left( A\right) >0\ ,\ \mathrm{tr}\left( A\right) >0.%
\end{array}%
\right.   \label{form_positive-2d}
\end{equation}
\end{quote}

\subsection{Dimension $n=3$}

Let $A$ be a self-adjoint operator in $\mathbb{C}^{3}$ and consider its
matrix representation with respect to the canonical basis (or with respect
to any other orthogonal basis):%
\begin{equation*}
\left[ A\right] =\left( a_{ij}\right) _{i,j=1,2,3}.
\end{equation*}%
As above, self-adjointness means%
\begin{equation*}
a_{ij}=\bar{a}_{ji}\in \mathbb{C}\ \  \forall i,j=1,2,3.
\end{equation*}%
The coefficients of the characteristic polynomial of $A$ are, as one can
obtain directly from (\ref{form_polinomio-characteristic}) or from Theorem %
\ref{teorema_coeficientes-pc}:%
\begin{eqnarray*}
b_{0} &=&1; \\
b_{1} &=&-\mathrm{tr}\left( A\right); \\
b_{2} &=&-\frac{1}{2}\left( b_{1}\mathrm{tr}\left( A\right) +b_{0}\mathrm{tr}%
\left( A^{2}\right) \right) =\frac{1}{2}\left( \mathrm{tr}\left( A\right)
^{2}-\mathrm{tr}\left( A^{2}\right) \right); \\
b_{3} &=&-\frac{1}{3}\left\{ b_{2}\mathrm{tr}\left( A\right) +b_{1}\mathrm{tr%
}\left( A^{2}\right) +b_{0}\mathrm{tr}\left( A^{3}\right) \right\} =-\det
\left( A\right).
\end{eqnarray*}%
The Hurwitz determinants of the characteristic polynomial are:%
\begin{equation*}
\Delta _{1}=b_{1}=-\mathrm{tr}\left( A\right),
\end{equation*}%
\begin{equation*}
\Delta _{2}=b_{1}b_{2}-b_{0}b_{3}=-\frac{1}{2}\mathrm{tr}\left( A\right)
\left( \mathrm{tr}\left( A\right) ^{2}-\mathrm{tr}\left( A^{2}\right)
\right) +\det \left( A\right),
\end{equation*}%
\begin{equation*}
\Delta _{3}=b_{1}b_{2}b_{3}-b_{0}b_{3}^{2}=\left\{ \frac{1}{2}\mathrm{tr}%
\left( A\right) \left( \mathrm{tr}\left( A\right) ^{2}-\mathrm{tr}\left(
A^{2}\right) \right) -\det \left( A\right) \right\} \det \left( A\right).
\end{equation*}%
Now, to get the necessary and sufficient conditions for $A$ be positive, we
have to consider four different situations, distinguished by the
multiplicity of the number zero as a root of the characteristic polynomial
of $A$. To save space (and also because it would be a little tedious), I
consider only the case in which the number zero is not a root of the
characteristic polynomial ($b_{3}=-\det \left( A\right) \neq 0$). In this
case, Theorem \ref{teorema_caracterizacao-op} stablishes that $A$ is
positive, if and only if,%
\begin{equation*}
\mathrm{tr}\left( A\right) >0,\ \ \frac{1}{2}\mathrm{tr}\left( A\right)
\left( \mathrm{tr}\left( A\right) ^{2}+\mathrm{tr}\left( A^{2}\right)
\right) +\det \left( A\right) >0 ,\ \ \left\{ \frac{1}{2}\mathrm{tr}%
\left( A\right) \left( \mathrm{tr}\left( A\right) ^{2}-\mathrm{tr}\left(
A^{2}\right) \right) -\det \left( A\right) \right\} \det \left( A\right) >0.
\end{equation*}

\begin{acknowledgement}
I thank my friend Wescley Bonomo for pointing out the Routh-Hurwitz
criterion.
\end{acknowledgement}

\end{document}